\documentclass[11pt,leqno]{amsart}
\usepackage{amscd}
\usepackage{amssymb}
\usepackage{amsfonts}
\usepackage{latexsym}
\usepackage{verbatim}

\theoremstyle{plain}
\newtheorem{theorem}{Theorem}[section]

\newtheorem{lemma}[theorem]{Lemma}
\newtheorem{prop}[theorem]{Proposition}
\newtheorem{cor}[theorem]{Corollary}
\newtheorem{rem}[theorem]{Remark}

\begin{document}
\title{On the Hermitian curvature of  symplectic manifolds}
\author{Luigi Vezzoni}
\date{\today}
\address{Dipartimento di Matematica ''L. Tonelli''\\ Universit\`a di Pisa\\
Largo Bruno Pontecorvo 5\\
56127 Pisa\\ Italy} \email{vezzoni@mail.dm.unipi.it}
\thanks{This work was supported by G.N.S.A.G.A.
of I.N.d.A.M.}
\begin{abstract}
In this paper we give conditions for the integrability of almost
complex structures
calibrated by symplectic forms.\\
We show that in the
symplectic case
Newlander-Nirenberg theorem reduces to
$\nabla''N_J=0$ and we give integrability  conditions in terms of the
curvature and the
Hermitian curvature of the induced metric.
\end{abstract}
\maketitle
%%%%%%%%%%%%%%%%% SIMBOLI MATEMATICI %%%%%%%%%%%%%%%%%%%%%%%%%%%%%%%%%%%%%%%%%%%%%%%%%%
\newcommand\C{{\mathbb C}}
\newcommand\R{{\mathbb R}}
\newcommand\Z{{\mathbb Z}}
\newcommand\T{{\mathbb T}}
\newcommand{\ov}[1]{\overline{#1}}
\newcommand{\gr}[1]{\mathbf{#1}}
\newcommand{\ka}[1]{\kappa{#1}}
\newcommand\de[2]{{\frac{\partial #1}{\partial #2}}}
\newcommand\J{{J\in\mathcal{C_{\kappa}}(M)}}
\newcommand\w[1]{{\widetilde{#1}}}
%%%%%%%%%%%%%%%%%%%%%%%%%%%%%Inizio%%%Lavoro%%%%%%%%%%%%%%%%%%%%%%%%%%%%%%%%%%%%%%%%%%%%%%%%%%%%%%%%%%%%%%%%%%%%%%%%%%%%%%%%%%%%%%%%%%%%%%%%%%%%%%%%%%%%%%%%%%%%%%%%%%%%%%%%%%%%%%%%%%%%%%%%%%%%%%%%%%%%%%%%%%%%%%%%%%%%%%%%%%%%%%%%%%%%%%%%%%%%%%%%%%%%%%%%%%%%%%%%%%%%%%%%%%%%%%%%%%%%%%%%%%%%%%%%%%%%%%%%%%%%%%%%%%%%%%%%%%
\section{Introduction}
The interplay between complex and symplectic structures has been recently
studied by many authors. Indeed, on any symplectic manifold
$(M,\kappa)$ there exists a $\kappa$-calibrated almost complex
structure $J$, so that $(M,g,J,\kappa)$ is an almost K\"ahler manifold.\\
In the context of almost K\"ahler geometry it is natural to study the
integrability of the almost complex structure.\\
In \cite{Gold} Goldberg proved that, if the curvature operator of an
almost K\"ahler manifold $(M,g,J,\kappa)$ commutes with $J$, then
$(M,g,J,\kappa)$ is a K\"ahler manifold.\\
He also conjectured that an Einstein almost K\"ahler metric on a
compact manifold is a K\"ahler metric.\\
If the scalar curvature is nonnegative the conjecture has been proved
by Sekigawa in \cite{Sek}.\\
For a survey, other references and results on this topic we refer to
\cite{L}.

In this paper we show that, if $(M,g,J,\kappa)$ is an almost K\"ahler
manifold satisfying certain properties, then it is K\"ahler.
Namely we give some conditions on the derivative of the
Nijenhuis tensor and on the curvature respectively, in order that
$J$ is integrable.\\
In section 2 we start by recalling some facts and fixing some
notations.\\
In section 3 we prove that if the $(0,1)$-part of the covariant
derivative of the Nijenhuis tensor of $J$ vanishes, then $J$ is
integrable (theorem \ref{nablaN_J}). This result generalizes theorem 2 of
\cite{vezio}.\\
In  section 4 we consider three types of curvature tensors on
$(M,g,J,\kappa)$. Namely, the Riemann curvature $R$, the
curvature $\w{R}$ of the Hermitian connection $\w{\nabla}$ and the
tensor $R^J$ defined by
$$
R^{J}(X,Y,Z,W)=g(\nabla_XJ\nabla_YZ-\nabla_YJ\nabla_XZ-\nabla_{[X,Y]}JZ,JW)\,.
$$
We show that  if $R^J$ and $R$ have the same components along
certain directions,
then $J$ is integrable (theorem \ref{th1}).\\
Finally we prove that if the bisectional curvature of $g$ and the Hermitian
bisectional curvature coincides, then $g$ is a K\"ahler metric
(theorem \ref{th3}).\\
A key tool in the proof of our results is the existence of
generalized normal holomorphic frames (see \cite{vezio}).
\section{Preliminaries}
Let  $M$ be a $2n$-dimensional (real)  manifold.\\
A  \emph{symplectic structure} on $M$
is a closed non-degenerate 2-form
$\kappa$, i.e. $d\kappa=0$ and $\kappa^n\neq0$. The pair $(M,\kappa)$
is said to be  a \emph{symplectic manifold}.

An \emph{almost complex structure} on $M$ is  a smooth section $J$
of $End(TM)$, such that
$J^2=-Id$.\\
An almost complex structure $J$ is said to be \emph{integrable} if the
Nijenhuis tensor
$$
N_J(X,Y)=[JX,JY]-J[JX,Y]-J[X,JY]-[X,Y]
$$
vanishes. In view of the celebrated
Newlander-Nirenberg theorem,
$J$ is integrable if and only if it is induced by
a holomorphic structure.

An almost complex structure $J$ on $M$
induces a natural splitting of the complexified of the tangent bundle.
Indeed let $T_J^{1,0}M$, $T_J^{0,1}M$, be the eigenspaces relatively to
$i$ and $-i$ respectively; then
$TM\otimes\C=T_J^{1,0}M\oplus T_J^{0,1}M$.
The sections of $T^{1,0}_JM$, $T_J^{0,1}M$ are called vector fields
of type $(1,0)$
and $(0,1)$ respectively.\\
We have  that $\overline{T_J^{1,0}M}=T^{0,1}_JM$ and that the
map $X\mapsto X-iJX$ defines an isomorphism between $TM$ and
$T^{1,0}_JM$.

An almost complex structure $J$ is said to be $\kappa$-\emph{calibrated} if
$$
g[x](\cdot,\cdot)=\kappa[x](\cdot,J_x\cdot)
$$
is a Hermitian metric on $M$. In this case the triple $(g,J,\kappa)$
is said to be an \emph{almost K\"ahler structure} on $M$ and $(M,g,J,\kappa)$
an \emph{almost K\"ahler manifold} .\\
Let  us denote by $\mathcal{C}_{\kappa}(M)$ the set of the almost
complex structure on $M$ calibrated by $\kappa$.
It is known that $\mathcal{C}_{\kappa}(M)$ is a non-empty and
contractible set (see e.g. \cite{A}). Therefore any symplectic
manifold admits almost K\"ahler structures.

Let $(M,\kappa)$  be a symplectic manifold and
$J\in\mathcal{C}_{\kappa}(M)$:
if $J$ is integrable the triple $(g,J,\kappa)$ is said to be a
\emph{K\"ahler structure} on $M$ and $(M,g,J,\kappa)$
a \emph{K\"ahler
manifold} (see e.g. \cite{F}).

%%%%%%%%%%%%%%%%%%%%%%%%%%%%%%%%%%%%%%%%%%%%%%%%%%%%%%%%%%%%%%%%%%%%%%%%%%%%%%%%%%%%%%%%%%%%%%%%%%%%%%%%%%%%%%%%%%%%%%%%%%%%%%%%%%%%%%%%%%%%%%%%%%%%%%%%%%%%%%%%%%%%%%%%%%%%%%%%%%%%%%%%%%%%%%%%%%%%%%%%%%%%%%%%%%%%%%%%%%%%%%%%%%%%%%%%%%%%%%%%%%%%%%%%%%%%%%%%%%%%%%%%%%%%%%%%%%%%%%%%%%%%%%%%%%%%%%%%%%%%%%%%%%%%%%%%%%%%%%%%%%%%%%%%%%%%%%%%%%%%%%%%%%%%%%%%%%%%%%%%%%%%%%%%%%%%%%%%%%%%%%%%%%%%%%%%%%%%%%%
\section{Some integrability conditions}
In this section we give conditions on the covariant derivative of the
Nijenhuis tensor in order that the complex structure is integrable.

Let
$(M,\kappa,J,\kappa)$ be a $2n$-dimensional almost K\"ahler manifold.
We denote  by $\nabla$ the
Levi-Civita connection relatively to $g$ and by $R$
the curvature tensor of $g$.\\
\newline
We have the following (see \cite{J}  and \cite{vezio})
\begin{theorem}[Generalized normal holomorphic frames]
\label{gnhf}
For any point $o\in M$ there exists a local complex $(1,0)$-frame,
$\{Z_{1},\dots ,Z_{n}\}$ around $o$, satisfying the following conditions:
\begin{enumerate}
\item[1.] $\nabla_{k}{Z}_{\ov{i}}(o)=0$, $1\leq k,i \leq n$\,;
\vspace{0.3 cm}
\item[2.] $\nabla_{k}Z_{i}(o)$ is of type $(0,1)$, $1\leq k,i \leq n$\,;
\vspace{0.3 cm}
\item[3.] if $\;G_{r\ov{s}}:=g(Z_{r},{Z}_{\ov{s}})$, then:
      $G_{r\ov{s}}(o)=\delta_{rs}$, $d G_{r\ov{s}}[o]=0$\,;
\vspace{0.3 cm}
\item[4.] $\nabla_{r}\nabla_{\ov{k}}Z_{i}(o)=0$\,,
$1\leq r,k,i\leq n$\,,
\end{enumerate}
where  $\ov{Z}_i:=Z_{\ov{i}}$ and
$\nabla_{Z_i}Z_j:=\nabla_i Z_j$.
\end{theorem}
By definition
$\{Z_{1},\dots ,Z_{n}\}$ is said a \emph{generalized normal holomorphic
frame around $o$}.

We recall that the following fundamental relation holds:
\begin{equation}
\label{fondamental}
2g((\nabla_{X}J)Y,Z)=g(N_{J}(Y,Z),JX)\,.
\end{equation}
From the complex extension of (\ref{fondamental}) it  can be easily
proved the following (see \cite{vezio})
\begin{cor}
\label{corart1}
The covariant derivative along  $(0,1)$-vector fields
of the almost complex structure $J$ is a tensor of complex type
$(0,1)$, i.e.
\begin{equation*}
(\nabla_{Z_{\ov{i}}}J)Z_{j}=0
\end{equation*}
for any $Z_i,Z_j$ complex vector fields of type $(1,0)$.
\end{cor}
\vspace{0.3cm}
Denote by $\nabla''$ the (0,1)-part of the covariant derivative
associated to the  Levi-Civita
connection, i.e.
$$
\nabla''_{W}X=\nabla_{W^{0,1}}X\,,
$$
where $W^{0,1}$ denotes the natural  projection of $W$ on $T^{0,1}_JM$ .
As a first application of the generalized normal holomorphic
frames we have the following
\begin{theorem}
\label{nablaN_J}
Assume that
$$
\nabla''N_J=0\,;
$$
then $(M,g,J,\kappa)$ is a K\"ahler manifold.
\end{theorem}
\begin{proof}
Let $o$ be an arbitrary point in $M$ and let $\{Z_1,\dots,Z_n\}$ be
a generalized normal holomorphic frame around $o$.\\
Since $N_J(T^{1,0}_JM\times T^{1,0}_JM)\subset T^{0,1}_JM$
we have
$$
Z_{\ov{i}}g(N_{J}(Z_i,Z_k),Z_{\ov{k}})=0\,.
$$
Therefore, by the assumption $\nabla''N_J=0$, we get
$$
g(N_{J}(Z_i,Z_k),\nabla_{\ov{i}}Z_{\ov{k}})=0
$$
and consequently
\begin{equation}
\label{epnNJ1}
g(N_{J}(Z_i,Z_k),[Z_{\ov{i}},Z_{\ov{k}}])=0\,.
\end{equation}
Moreover a direct computation gives
\begin{equation}
\label{epnNJ2}
N_{J}(Z_i,Z_k)(o)=-\frac{1}{4}[Z_i,Z_k](o)\,.
\end{equation}
Hence (\ref{epnNJ1}) and (\ref{epnNJ2}) imply $N_{J}(Z_i,Z_k)(o)=0\,$,
i.e. $J$ is integrable, so that $(M,g,J,\kappa)$ is a K\"ahler manifold.
\end{proof}
\begin{rem}{\rm
Let $B$ the (2,1)-tensor on $(M,g,J,\kappa)$ defined by
\begin{equation}
\label{B}
B(X,Y)=J(\nabla_XJ)Y-(\nabla_{JX}J)Y\,.
\end{equation}

Then the Nijenhuis tensor of $J$ is the antisymmetric part of $B$.\\
In \cite{vezio} (theorem 2) it is proved that an almost complex structure $J$
is integrable if and only if $\nabla''B=0$ . Therefore theorem 3.2
generalizes theorem 2 of \cite{vezio} .}
\end{rem}
\vspace{0.3cm}
\begin{rem}{\rm
With  respect to a generalized normal holomorphic frame,
$\{Z_1,\dots,Z_n\}$, the components of the
curvature tensor  $R_{i\ov{k}r\ov{s}}(o)$ and
$R_{i\ov{k}\ov{r}\ov{s}}(o)$ are
given by
\begin{itemize}
\item
  $R_{i\ov{k}r\ov{s}}(o)=-g(\nabla_{\ov{k}}\nabla_iZ_r,Z_{\ov{s}})(o)$
  ,
\vspace{0.2 cm}
\item
  $R_{i\ov{k}\ov{r}\ov{s}}(o)=g(\nabla_{i}\nabla_{\ov{k}}Z_{\ov{r}},
Z_{\ov{s}})(o)$ .
\end{itemize}
}
\end{rem}
\vspace{0.3cm}
In order to study the integrability of a $\kappa$-calibrated almost
complex
structure $J$ it is useful to introduce the following tensor
$$
L(X,Y,Z,W)=g((\nabla_XB)(Y,Z),W)\,,
$$
where $B$ is defined by (\ref{B}).\\
A first property of the tensor $L$ is given by the following
\begin{lemma}
\label{L.1}
This identities hold
\begin{equation*}
\begin{aligned}
&L_{\ov{i}jr\ov{s}}=L_{j\ov{i}\ov{s}r}\,;\\
&L_{i{\ov{j}}rs}=L_{ij\ov{r}s}=
L_{ijr\ov{s}}=0\,,
\end{aligned}
\end{equation*}
for $1\leq i,j,r,s\leq n$.
\end{lemma}
\begin{proof}
Let $o$ be  an arbitrary point in $M$ and let
$\{Z_1,\dots,Z_n\}$ be a generalized normal holomorphic frame around $o$.
By the definition of $L$ and the properties of generalized normal
holomorphic frames, we get
\begin{equation*}
\begin{aligned}
L_{\ov{i}jr\ov{s}}(o)=&g((\nabla_{\ov{{i}}}B)(Z_j,Z_r),Z_{\ov{s}})(o)\\
                     =&g(\nabla_{\ov{{i}}}(B(Z_j,Z_r)),Z_{\ov{s}})(o)\\
=&g(\nabla_{\ov{i}}(J\nabla_jJ)Z_r,Z_{\ov{s}})(o)
-ig(\nabla_{\ov{i}}(\nabla_{j}J)Z_r,Z_{\ov{s}})(o)\\
=&2ig(\nabla_{\ov{i}}J\nabla_jZ_r,Z_{\ov{s}})(o)
+2g(\nabla_{\ov{i}}\nabla_jZ_r,Z_{\ov{s}})(o)\\
=&2ig(\nabla_{\ov{i}}J\nabla_jZ_r,Z_{\ov{s}})(o)+2R_{\ov{i}jr\ov{s}}(o)\,,
\end{aligned}
\end{equation*}
i.e.
\begin{equation}
\begin{aligned}
\label{FL.1}
L_{\ov{i}jr\ov{s}}(o)
=&2ig(\nabla_{\ov{i}}J\nabla_jZ_r,Z_{\ov{s}})(o)+2R_{\ov{i}jr\ov{s}}(o)\,.
\end{aligned}
\end{equation}
By a direct computation we get
\begin{equation*}
\begin{aligned}
L_{\ov{i}jr\ov{s}}(o)=&
2i Z_{\ov{i}}g(J\nabla_jZ_r,\ov{Z}_s)(o)
+2ig(\nabla_jZ_r,J\nabla_{\ov{i}}Z_{\ov{s}})(o)\\
&+2R_{\ov{i}jr\ov{s}}(o)\,.
\end{aligned}
\end{equation*}
Then we have
\begin{equation*}
\begin{aligned}
L_{\ov{i}jr\ov{s}}(o)=&
2i Z_{\ov{i}}g(J\nabla_jZ_r,Z_{\ov{s}})(o)
+2iZ_jg(Z_r,J\nabla_{\ov{i}}Z_{\ov{s}})(o)\\
&-2ig(Z_r,\nabla_jJ\nabla_{\ov{i}}Z_{\ov{s}})(o)+2R_{\ov{i}jr\ov{s}}(o)\\
=&-2Z_{\ov{i}}g(\nabla_jZ_r,Z_{\ov{s}})(o)+
2Z_jg(Z_r,\nabla_{\ov{i}}Z_{\ov{s}})(o)\\
&-2ig(Z_r,\nabla_jJ\nabla_{\ov{i}}Z_{\ov{s}})(o)+2R_{\ov{i}jr\ov{s}}(o)\\
=&-2g(\nabla_{\ov{i}}\nabla_jZ_r,Z_{\ov{s}})(o)+2g(Z_r,\nabla_j\nabla_{\ov{i}}Z_{\ov{s}})(o)\\
 &-2ig(Z_r,\nabla_jJ\nabla_{\ov{i}}Z_{\ov{s}})(o)+2 R_{j\ov{i}\ov{s}r}(o)\\
=&-2R_{\ov{i}jr\ov{s}}(o)+2R_{j\ov{i}\ov{s}r}(o)\\
 &-2ig(Z_r,\nabla_jJ\nabla_{\ov{i}}Z_{\ov{s}})(o)+2 R_{j\ov{i}\ov{s}r}(o)\\
=&-2ig(Z_r,\nabla_jJ\nabla_{\ov{i}}Z_{\ov{s}})(o)+2 R_{j\ov{i}\ov{s}r}(o)\,.
\end{aligned}
\end{equation*}
In a similar way we obtain
\begin{equation*}
L_{j\ov{i}\ov{s}r}(o)=
-2ig(Z_r,\nabla_jJ\nabla_{\ov{i}}Z_{\ov{s}})(o)+2 R_{j\ov{i}\ov{s}r}(o)\,.
\end{equation*}
Therefore
\begin{equation*}
L_{\ov{i}jr\ov{s}}(o)=L_{j\ov{i}\ov{s}r}(o)\,.
\end{equation*}
The second part of the proof is straightforward.
\end{proof}
As a  corollary of theorem \ref{nablaN_J}
we have the following (see also \cite{vezio}):
\begin{prop}
Assume that
$$
L(Z_{\ov{i}},Z_i,Z_j,Z_{\ov{j}})=0\,,
$$
for any complex vector fields $Z_i,Z_j$ of type $(1,0)$;
then $(M,g,J,\kappa)$ is a K\"ahler manifold.
\end{prop}
\begin{proof}
Let $\{Z_1,\dots,Z_n\}$ be a generalized normal holomorphic frame
around $o$. A direct  computation gives
\begin{equation}
\label{L=0}
L(Z_{\ov{i}},Z_i,Z_j,Z_{\ov{j}})(o)=
g(\nabla_{i}Z_j,\nabla_{\ov{i}}Z_{\ov{j}})(o)\,;
\end{equation}
hence $\nabla_iZ_j(o)=0$ .\\
Therefore $N_J(Z_i,Z_j)(o)=-\frac{1}{4}[Z_i,Z_j](o)=0$ .
\end{proof}
\section{Curvature and integrability}
In this section we give an integrability condition in terms of curvature.

We start defining the following tensor
$$
R^{J}(X,Y,Z,W)=g(\nabla_XJ\nabla_YZ-\nabla_YJ\nabla_XZ-\nabla_{[X,Y]}JZ,JW)\,.
$$
We have
\begin{lemma}
Let $\{Z_1,\dots,Z_n\}$ be  a generalized normal holomorphic frame around
$o$, then
$$
R^J_{\ov{i}jr\ov{s}}(o)=-ig(\nabla_{\ov{i}}J\nabla_jZ_r,Z_{\ov{s}})(o)\,.
$$
\end{lemma}
\begin{proof}
By the definition  of $R^J$ it follows that
\begin{equation*}
\begin{aligned}
R^J_{\ov{i}jr\ov{s}}(o)=&-ig(\nabla_{\ov{i}}J\nabla_jZ_r,Z_{\ov{s}})(o)+ig(\nabla_{j}J\nabla_{\ov{i}}Z_r,Z_{\ov{s}})(o)\\
                        &-g(\nabla_{[Z_{\ov{i}},Z_j]}Z_r,Z_{\ov{s}})(o)\\
                       =&-ig(\nabla_{\ov{i}}J\nabla_jZ_r,Z_{\ov{s}})(o)+ig(\nabla_{j}J\nabla_{\ov{i}}Z_r,Z_{\ov{s}})(o)\,.
\end{aligned}
\end{equation*}
By corollary \ref{corart1}, we have
$(\nabla_{\ov{i}}J)Z_j=0$ . Therefore  we obtain
$$
g(\nabla_{j}J\nabla_{\ov{i}}Z_r,Z_{\ov{s}})(o)=ig(\nabla_{j}\nabla_{\ov{i}}Z_r,Z_{\ov{s}})(o)=0\,
$$
and then
$$
R^J_{\ov{i}jr\ov{s}}(o)=-ig(\nabla_{\ov{i}}J\nabla_jZ_r,Z_{\ov{s}})(o)\,.
$$
\end{proof}
The previous lemma and equation (\ref{FL.1}) give us
\begin{equation}
\label{FL.2}
L_{\ov{i}jr\ov{s}}(o)=-2R^{J}_{\ov{i}jr\ov{s}}(o)+2R_{\ov{i}jr\ov{s}}(o)\,.
\end{equation}
Moreover  equation (\ref{FL.2}) and lemma (\ref{L.1}) imply the
following
\begin{theorem}
\label{th1}
If
$$
R^{J}(Z_{\ov{i}},Z_j,Z_i,Z_{\ov{i}})=R(Z_{\ov{i}},Z_j,Z_i,Z_{\ov{j}})\,,
$$
for any $Z_i,Z_j$ of type $(1,0)$, then $J$ is integrable.
\end{theorem}
\vspace{0.3 cm}
\begin{rem}
\rm{The hypothesis of the previous theorem can be replaced by
$$
R^{J}(JX,Y,X,JY)=R(JX,Y,X,JY)
$$
for any  $X,Y$ real vector fields.}
\end{rem}
\vspace{0.3 cm}
Let $\w{\nabla}$ be the
\emph{Hermitian connection} on $(M,g,J,\kappa)$; $\w{\nabla}$ is the
connection defined by
$$
\w{\nabla}=\nabla-\frac{1}{2}J\nabla J\,.
$$
Then $\w{\nabla}$ preserves the metric $g$, the almost complex
structure $J$ and its torsion is given by the Nijenhuis torsion,
namely
$$
\w{\nabla}g=0,\quad \w{\nabla}J=0\,,\quad T^{\w{\nabla}}
=\frac{1}{4}N_J\,.
$$
Let us denote by $\w{R}$ the curvature tensor
of $\w{\nabla}$.\\
From \cite{Sek} we have the following
\begin{lemma}[\cite{Sek}]
\label{lemmasek}
For any $X,Y,Z,W$ we have the following formula
\begin{equation*}
\begin{aligned}
\w{R}(X,Y,Z,W)=&\frac{1}{2}R(X,Y,Z,W)+\frac{1}{2}R(X,Y,JZ,JW)\\
&-\frac{1}{4}g((\nabla_XJ)(\nabla_YJ)Z-(\nabla_YJ)(\nabla_XJ)Z,W)\,.
\end{aligned}
\end{equation*}
\end{lemma}
Now we have
\begin{lemma}
\label{wR1}
Let $\{Z_1,\dots,Z_n\}$ be a local $(1,0)$-frame. Then
\begin{equation}
\label{corR}
\w{R}_{i\ov{j}r\ov{s}}=R_{i\ov{j}r\ov{s}}+\frac{1}{4}L_{i\ov{j}r\ov{s}}
\end{equation}
for any $1\leq i,j,r,s\leq n$.
\end{lemma}
\begin{proof}
Let $\{Z_1,\dots,Z_n\}$ be  a local $(1,0)$-frame; from
lemma \ref{lemmasek} we have
\begin{equation*}
\label{corR}
\begin{aligned}
\w{R}_{i\ov{j}r\ov{s}}=&\frac{1}{2}R_{i\ov{j}r\ov{s}}+
\frac{1}{2}R_{i\ov{j}r\ov{s}}
-\frac{1}{4}g((\nabla_{i}J)(\nabla_{\ov{j}}J)Z_r,Z_{\ov{s}})
+\frac{1}{4}g((\nabla_{\ov{j}}J)(\nabla_iJ)Z_{r},Z_{\ov{s}})
\\
=&R_{i\ov{j}r\ov{s}}+
\frac{1}{4}g((\nabla_{\ov{j}}J)(\nabla_iJ)Z_{r},Z_{\ov{s}})\,.
\end{aligned}
\end{equation*}
From corollary \ref{corart1} we have  $(\nabla_{\ov{j}}J)Z_r$;
then we get
$$
g((\nabla_{i}J)(\nabla_{\ov{j}}J)Z_r,Z_{\ov{s}})=0\,.
$$
Therefore
\begin{equation*}
\begin{aligned}
\w{R}_{i\ov{j}r\ov{s}}
=&R_{i\ov{j}r\ov{s}}+
i\frac{1}{4}g((\nabla_{\ov{j}}J)\nabla_iZ_{r},Z_{\ov{s}})
-\frac{1}{4}g((\nabla_{\ov{j}}J)J\nabla_iZ_{r},Z_{\ov{s}})\\
=&R_{i\ov{j}r\ov{s}}+
i\frac{1}{4}g(\nabla_{\ov{j}}J\nabla_iZ_{r},Z_{\ov{s}})
-i\frac{1}{4}g(J\nabla_{\ov{j}}\nabla_iZ_{r},Z_{\ov{s}})
+\frac{1}{4}g((\nabla_{\ov{j}}\nabla_iZ_{r},Z_{\ov{s}})\\
&+\frac{1}{4}g((J\nabla_{\ov{j}}J\nabla_iZ_{r},Z_{\ov{s}})\\
=&R_{i\ov{j}r\ov{s}}+
i\frac{1}{4}g(\nabla_{\ov{j}}J\nabla_iZ_{r},Z_{\ov{s}})
+\frac{1}{4}g(\nabla_{\ov{j}}\nabla_iZ_{r},Z_{\ov{s}})
+\frac{1}{4}g((\nabla_{\ov{j}}\nabla_iZ_{r},Z_{\ov{s}})\\
&+i\frac{1}{4}g((\nabla_{\ov{j}}J\nabla_iZ_{r},Z_{\ov{s}})\\
=&R_{i\ov{j}r\ov{s}}+
i\frac{1}{2}g(\nabla_{\ov{j}}J\nabla_iZ_{r},Z_{\ov{s}})
+\frac{1}{2}g(\nabla_{\ov{j}}\nabla_iZ_{r},Z_{\ov{s}})\,,
\end{aligned}
\end{equation*}
i.e.
\begin{equation}
\label{wR}
\w{R}_{i\ov{j}r\ov{s}}=R_{i\ov{j}r\ov{s}}+
i\frac{1}{2}g(\nabla_{\ov{j}}J\nabla_iZ_{r},Z_{\ov{s}})
+\frac{1}{2}R_{\ov{j}ir\ov{s}}\,.
\end{equation}
Let assume now that $\{Z_1,\dots,Z_n\}$ is a generalized normal
holomorphic frame around a point $o$. Then equations  (\ref{wR}) and
(\ref{FL.1}) imply
\begin{equation*}
\begin{aligned}
\w{R}_{i\ov{j}r\ov{s}}(o)=R_{i\ov{j}r\ov{s}}(o)
+\frac{1}{4}L_{\ov{j}ir\ov{s}}(o)\,.
\end{aligned}
\end{equation*}
\end{proof}
Theorem \ref{th1} and lemma \ref{wR1} imply the following
\begin{theorem}
\label{th2}
If
$$
\w{R}(Z_i,Z_{\ov{i}},Z_j,Z_{\ov{j}})=R(Z_i,Z_{\ov{i}},Z_j,Z_{\ov{j}})
$$
for any $Z_i,Z_j$ $(1,0)$-fields, then $(M,g,J,\kappa)$ is a  K\"ahler
manifold.
\end{theorem}
\vspace{0.3 cm}
\begin{rem}
\rm{The hypothesis of theorem \ref{th2} can be given as
$$
\w{R}(X,JX,Y,JY)=R(X,JX,Y,JY)\,.
$$
for any pair of real vector fields $X,Y$}
\end{rem}
\vspace{0.3 cm}
The last theorem can be stated in terms of the bisectional curvatures
of $\nabla$ and $\w{\nabla}$.\\
Let $p\in M$ and $v,w$ be two unit vector of $T_pM$. Then the
holomorphic bisectional curvature of the planes $\sigma_1=<v,Jv>$,
$\sigma_2=<w,Jw>$  is defined by
$$
K[p](\sigma_1,\sigma_2)=R(v,Jv,w,Jw)\,.
$$
We denote by $\w{K}[p](\sigma_1,\sigma_2)$ the Hermitian bisectional
curvature (i.e. the bisectional curvature of the Hermitian connection).
For any $p\in M$  let $\mathcal{P}^{1,1}_p(M)$ be the set of $J$-invariant
planes in $T_pM$.\\
Hence we get
\begin{theorem}
\label{th3}
If
$$
K(p)(\sigma_1,\sigma_2)=\w{K}(p)(\sigma_1,\sigma_2)
$$
for any $p\in M$ and $\sigma_1,\sigma_2\in\mathcal{P}^{1,1}_p(M)$,
then $J$ is an integrable almost complex structure and therefore
$(M,g,J,\kappa)$ is a K\"ahler manifold.
\end{theorem}
\vspace{0.3 cm}
%\textbf{Acknowledgments:} The author

\end{document}